# Decoupled, Energy Stable Scheme for Hydrodynamic Allen-Cahn Phase Field Moving Contact Line Model[*]

Rui Chen

*Institute of Applied Physics and Computational Mathematics, Beijing 100088, P.R. China*
*School of Science, Beijing University of Posts and Telecommunications, Beijing 100876, P.R. China*
*Email: ruichenbnu@gmail.com*

Xiaofeng Yang

*Department of Mathematics, University of South Carolina, Columbia, SC, 29208, USA*
*Email: xfyang@math.sc.edu*

Hui Zhang

*School of Mathematical Sciences, Beijing Normal University, Laboratory of Mathematics and Complex Systems, Ministry of Education, Beijing 100875, P.R. China*
*Email: hzhang@bnu.edu.cn*

**Abstract**

In this paper, we present an efficient energy stable scheme to solve a phase field model incorporating contact line condition. Instead of the usually used Cahn-Hilliard type phase equation, we adopt the Allen-Cahn type phase field model with the static contact line boundary condition that coupled with incompressible Navier-Stokes equations with Navier boundary condition. The projection method is used to deal with the Navier-Stokes equations and an auxiliary function is introduced for the non-convex Ginzburg-Landau bulk potential. We show that the scheme is linear, decoupled and energy stable. Moreover, we prove that fully discrete scheme is also energy stable. An efficient finite element spatial discretization method is implemented to verify the accuracy and efficiency of proposed schemes. Numerical results show that the proposed scheme is very efficient and accurate.

## 1. Introduction

Two phase immiscible flows are common in our lives, such as the air bubble in the water, and droplet of oil in the water, and cells in the blood, etc.. When the air bubble rises to the edge of glass, the interface will change its geometry. Then it will form a moving contact line (MCL) problem. Moving contact line, where the fluid-fluid interface touches the solid wall, has been widely investigated by the researcher theoretically and experimentally. People are interested in the topic that how the contact line evolves when the solid wall moves. In this situation, the no-slip boundary condition for the flow is not applicable again. From this viewpoint, Qian et al. [28, 29] proposed the generalized Navier boundary condition (GNBC) according to the molecular dynamics theory. To investigate the complex behavior at MCLs, plenty of models are studied by the researchers. For example, molecular dynamical (MD) simulations are studied by Koplik et al. [15, 16] and Qian et al. [28, 29]. Microscopic-macroscopic hybrid simulations had been carried out by Hadjiconstantinou [8], Ren and E [30] etc. Although this approach is powerful, the computational cost is very expensive for macroscopic applications.

There are a lot of methods to study the MCL problem, such as the immersed interface method [17], the volume of fluid method [31], a hybrid atomistic-continuum method [8], and phase field method [41] etc. In the recent years, phase field method has been used widely on the interfacial phenomena, and has been applied to simulate many dynamical processes successfully in many fields [2, 4, 11–14, 18–20, 22, 25, 27, 35, 36, 38–43, 46–49, 51, 52, 56, 61, 63, 64]. In a phase field model, a continuous phase field function is used to denote the two immiscible fluids where the fluid-fluid interface has a thickness. In the framework of phase field approach, the governing system is usually obtained from the gradient flow, which is a variational formalism from the total free energy. Thus, there is a physical energy law associated with the phase field model [29]. This energy law will help us to carry out mathematical analysis and further design efficient numerical scheme. Thus, from the numerical point view, we are interested in constructing a simple and efficient (linear and decoupled) scheme which satisfies the discrete energy law.

Qian et al. [29] have proposed a Navier-Stokes Cahn-Hilliard (NSCH) system to study the MCL problem, where the GNBC is used for Navier-Stokes equation and the dynamic contact line condition (DCLC) is applied on the Cahn-Hilliard equation. As we all know, the Cahn-Hilliard equation is a four-order equation with volume fraction being conserved, but it needs more time to carry out than Allen-Cahn (second-order) equation on numerical computation. Though Allen-Cahn equation is not conserved on volume fraction, it can keep conserved with the introduction of a Lagrangian multiplier. Thus in this paper we use Allen-Cahn equation to replace the Cahn-Hilliard equation.

The purpose of our paper is to construct a linear, decoupled, fully discrete, first-order, and unconditionally energy stable scheme for the Navier-Stokes Allen-Cahn (NSAC) system. We know that the NSAC system is a coupled system, which includes the coupling between the phase field variable and the velocity in the convection and the stress, and the coupling between the velocity and the pressure in the momentum equation. There are several methods to design the energy stable scheme, such as the operator-splitting on the time-discretization [10], a convex splitting scheme in [1, 33], and the stabilization approach in [21, 22, 24, 36–39]. Although these schemes are energy stable, they still have some disadvantages. Firstly, the operator-splitting scheme and the convex splitting scheme are usually coupled and nonlinear due to the couplings in the system, which take much time on iterations to carry out the numerical results. Secondly, it is difficult to prove the unconditional solvability for these nonlinear schemes. At last, although the stabilizing approach is linear and decoupled, the truncation error is introduced in the scheme, which requires the double well potential to be bounded form the viewpoint of mathematics. Thus, in this paper, we shall overcome these difficulties to construct a fully discrete, linear and decoupled energy stable scheme using the recently developed "Invariant Energy Quadratization" approach [5, 9, 50, 53–55, 57, 58, 60, 62]. In our previous paper [24, 59], we developed some decoupled, energy stable numerical schemes to solve the hydrodynamics coupled Allen-Cahn MCL phase field model by adding the extra stabilizing terms for phase field equation. In this paper we adopt a novel skill to construct a linear, decoupled, and unconditionally energy stable scheme without using the stabilizing approach for the phase field equation.

The rest of the paper is organized as follows. In the next section, we present the phase field model of moving contact line condition and show the energy dissipation law for the system. In section 3, we construct a linear decoupled energy stable scheme for this coupled nonlinear system by the introduction of the auxiliary function. In section 4, we implement the finite element method on this model for spatial discretization and prove that the fully discrete scheme is energy stable. In section 5, we present some numerical simulations to illustrate the efficiency and accuracy of the proposed numerical scheme. Finally, we give some conclusions in the last section.

## 2. A Navier-Stokes Allen-Cahn coupled model

In [28, 29], the Navier-Stokes Cahn-Hilliard coupled system (NSCH) with a GNBC was developed to study the two-phase incompressible, immiscible fluid, where the fluid-fluid interface touches the solid wall. A non-dimensional version of the system is given as follows.

Incompressible Navier-Stokes equations for hydrodynamics

$$\mathbf{u}_t + (\mathbf{u} \cdot \nabla)\mathbf{u} = \nu\Delta\mathbf{u} - \nabla p + \lambda\mu\nabla\phi, \tag{2.1}$$

$$\nabla \cdot \mathbf{u} = 0, \tag{2.2}$$

$$\mathbf{u} \cdot \mathbf{n} = 0, \quad \text{on} \quad \partial\Omega, \tag{2.3}$$

$$l(\phi)(\mathbf{u}_\tau - \mathbf{u}_w) + \nu\partial_{\mathbf{n}}\mathbf{u}_\tau - \lambda L(\phi)\nabla_\tau\phi = 0, \quad \text{on} \quad \partial\Omega. \tag{2.4}$$

Cahn-Hilliard type phase field equations:

$$\phi_t + \mathbf{u} \cdot \nabla\phi = M\Delta\mu, \tag{2.5}$$

$$\mu = -\varepsilon\Delta\phi + f(\phi), \tag{2.6}$$

$$\partial_{\mathbf{n}}\mu = 0, \quad \text{on} \quad \partial\Omega, \tag{2.7}$$

$$\phi_t + \mathbf{u}_\tau \cdot \nabla_\tau\phi = -\gamma L(\phi), \quad \text{on} \quad \partial\Omega, \tag{2.8}$$

where $\mathbf{u}$ is the fluid velocity, $p$ is the pressure, $\phi$ is the phase field variable, $\mu$ is the chemical potential, the function $L(\phi)$ is given by

$$L(\phi) = \varepsilon\partial_{\mathbf{n}}\phi + g'(\phi), \tag{2.9}$$

where $g(\phi)$ is the boundary interfacial energy, $l(\phi) \geq 0$ is a given coefficient function. The function $f(\phi) = F'(\phi)$ with $F(\phi)$ being the Ginzburg-Landau bulk potential. More precisely, $F(\phi)$ and $g(\phi)$ are defined as

$$F(\phi) = \frac{1}{4\varepsilon}(\phi^2 - 1)^2, \quad g(\phi) = -\frac{\sqrt{2}}{3}\cos\theta_s\sin(\frac{\pi}{2}\phi), \tag{2.10}$$

where $\theta_s$ is the static contact angle. In equations (2.1)-(2.8), $\nabla$ denote the gradient operator, $\mathbf{n}$ is the outward normal direction on boundary $\partial\Omega$, $\tau$ is tangential direction on the boundary, and vector operator $\nabla_\tau = \nabla - (\mathbf{n} \cdot \nabla)\mathbf{n}$ is the gradient along tangential direction, $\mathbf{u}_w$ is the boundary wall velocity, $\mathbf{u}_\tau$ is the boundary fluid velocity in tangential direction. From (2.3), we have $\mathbf{u} = \mathbf{u}_\tau$ on boundary $\partial\Omega$. There are six non-dimensional parameters in this system. $\nu$ is the viscosity coefficient, $\lambda$ denotes the strength of the capillary force comparing to the Newtonian fluid stress, $M$ is the mobility coefficient, $\gamma$ is a boundary relaxation coefficient, $l(\phi)$ is the ratio of domain size to boundary slip length, $\varepsilon$ denotes the interface thickness.

When $\gamma \to +\infty$, the DCLC (2.8) reduces to static contact line condition (SCLC),

$$L(\phi) = 0, \quad \text{on} \quad \partial\Omega, \tag{2.11}$$

and the GNBC reduces to the Navier boundary condition (NBC),

$$l(\phi)(\mathbf{u}_\tau - \mathbf{u}_w) + \nu \partial_\mathbf{n} \mathbf{u}_\tau = 0, \quad \text{on} \quad \partial\Omega. \tag{2.12}$$

If we further set $g'(\phi) \equiv 0$, the SCLC (2.11) turns to the Neumann boundary condition. If we take $l(\phi) \to +\infty$ in (2.12), then the Navier slip boundary condition reduces to traditional no-slip boundary condition.

In this paper, we replace the fourth order Cahn-Hilliard equation with the second order Allen-Cahn equation with a Lagrangian multiplier. In this new Navier-Stokes Allen-Cahn (NSAC) system, the fluid equation is same as (2.1)-(2.4). The non-dimensional, Allen-Cahn phase field equation is given as follows.

$$\phi_t + \mathbf{u} \cdot \nabla\phi = -M\mu, \tag{2.13}$$

$$\mu = -\varepsilon\Delta\phi + f(\phi), \tag{2.14}$$

$$\phi_t + \mathbf{u}_\tau \cdot \nabla_\tau \phi = -\gamma L(\phi), \quad \text{on} \quad \partial\Omega. \tag{2.15}$$

To preserve the volume fraction, we add the Lagrangian multiplier $\xi(t)$ as follows.

$$\mu = -\varepsilon\Delta\phi + f(\phi) + \xi(t), \tag{2.16}$$

$$\frac{d}{dt}\int_\Omega \phi dx = 0. \tag{2.17}$$

**Remark 2.1.** *By the integration of* (2.13) *on* $\Omega$, *we can obtain the Lagrangian multiplier* $\xi(t)$ *depending on* $t$,

$$\xi(t) = -\frac{1}{|\Omega|}\Big(\int_\Omega f(\phi)dx - \int_{\partial\Omega} \varepsilon\partial_\mathbf{n}\phi ds\Big). \tag{2.18}$$

We now derive the energy dissipation law for PDEs system (2.1)-(2.4) and (2.13), (2.15), (2.16). Here and after, for any function $f, g \in L^2(\Omega)$, we use $(f, g)$ to denote $\int_\Omega fg dx$, $(f, g)_{\partial\Omega}$ to denote $\int_{\partial\Omega} fg dx$, and $\|f\|^2 = (f, f)$ and $\|f\|^2_{\partial\Omega} = (f, f)_{\partial\Omega}$.

**Theorem 2.1.** *The NSAC system* (2.1)-(2.4) *and* (2.13), (2.15), (2.16) *with GNBC* (2.4) *and DCLC* (2.8) *is a dissipative system satisfying the following energy dissipation law,*

$$\frac{dE_{tot}}{dt} = -\nu\|\nabla\mathbf{u}\|^2 - \lambda M\|\mu\|^2 - \lambda\gamma\|L(\phi)\|^2_{\partial\Omega} - \|l(\phi)^{\frac{1}{2}}\mathbf{u}_s\|^2_{\partial\Omega} - (l(\phi)\mathbf{u}_s, \mathbf{u}_w)_{\partial\Omega}, \tag{2.19}$$

*where* $\mathbf{u}_s = \mathbf{u}_\tau - \mathbf{u}_w$ *is the velocity slip on boundary* $\partial\Omega$, $E_{tot} = E_k[\mathbf{u}] + E_b[\phi] + E_s[\phi]$, *and*

$$E_k[\mathbf{u}] = \frac{\|\mathbf{u}\|^2}{2}, \quad E_b[\phi] = \lambda\varepsilon\frac{\|\nabla\phi\|^2}{2} + \lambda(F(\phi), 1), \quad E_s[\phi] = \lambda(g(\phi), 1)_{\partial\Omega}. \tag{2.20}$$

*Proof.* By taking the inner product of (2.1) with $\mathbf{u}$, using the incompressible condition (2.2) and the zero flux boundary condition (2.3), we have

$$\frac{1}{2}\frac{d}{dt}\|\mathbf{u}\|^2 = \nu(\partial_\mathbf{n}\mathbf{u}, \mathbf{u})_{\partial\Omega} - \nu\|\nabla\mathbf{u}\|^2 + \lambda(\mu\nabla\phi, \mathbf{u}). \tag{2.21}$$

By taking the inner product of (2.13) with $\lambda\mu$, we get

$$\lambda(\phi_t, \mu) + \lambda(\mathbf{u} \cdot \nabla\phi, \mu) = -\lambda M\|\mu\|^2. \tag{2.22}$$

By taking the inner product of (2.16) with $\lambda\phi_t$, we have

$$\lambda(\mu, \phi_t) = -\lambda\varepsilon(\partial_{\mathbf{n}}\phi, \phi_t)_{\partial\Omega} + \frac{1}{2}\lambda\varepsilon\frac{d}{dt}\|\nabla\phi\|^2 + \lambda\frac{d}{dt}(F(\phi), 1). \tag{2.23}$$

Combining equations (2.21)-(2.23), we obtain

$$\begin{aligned}
\frac{1}{2}\frac{d}{dt}\|\mathbf{u}\|^2 + \frac{1}{2}\lambda\varepsilon\frac{d}{dt}\|\nabla\phi\|^2 + \lambda\frac{d}{dt}(F(\phi), 1) &= -\nu\|\nabla\mathbf{u}\|^2 - \lambda M\|\mu\|^2 \\
&\quad + \nu(\partial_{\mathbf{n}}\mathbf{u}, \mathbf{u})_{\partial\Omega} + \lambda\varepsilon(\partial_{\mathbf{n}}\phi, \phi_t)_{\partial\Omega}. \tag{2.24}
\end{aligned}$$

Then, by using (2.8), (2.9) and boundary condition (2.4), we have

$$\begin{aligned}
\nu(\partial_{\mathbf{n}}\mathbf{u}, \mathbf{u})_{\partial\Omega} = \nu(\partial_{\mathbf{n}}\mathbf{u}_\tau, \mathbf{u}_\tau)_{\partial\Omega} &= (\lambda L(\phi)\nabla_\tau\phi - l(\phi)(\mathbf{u}_\tau - \mathbf{u}_w), \mathbf{u}_\tau)_{\partial\Omega} \\
&= \lambda(L(\phi)\nabla_\tau\phi, \mathbf{u}_\tau)_{\partial\Omega} - (l(\phi)\mathbf{u}_s, \mathbf{u}_s + \mathbf{u}_w)_{\partial\Omega}, \tag{2.25}
\end{aligned}$$

and

$$\begin{aligned}
\lambda\varepsilon(\partial_{\mathbf{n}}\phi, \phi_t)_{\partial\Omega} &= \lambda(L(\phi) - g'(\phi), \phi_t)_{\partial\Omega} \\
&= \lambda(L(\phi), \phi_t)_{\partial\Omega} - \lambda(g'(\phi), \phi_t)_{\partial\Omega} \\
&= \lambda(L(\phi), -\mathbf{u}_\tau \cdot \nabla_\tau\phi - \gamma L(\phi))_{\partial\Omega} - \lambda\frac{d}{dt}(g(\phi), 1)_{\partial\Omega} \\
&= -\lambda(L(\phi)\nabla_\tau\phi, \mathbf{u}_\tau)_{\partial\Omega} - \lambda\gamma\|L(\phi)\|_{\partial\Omega}^2 - \lambda\frac{d}{dt}(g(\phi), 1)_{\partial\Omega}. \tag{2.26}
\end{aligned}$$

Summing up (2.24), (2.25) and (2.26), we obtain the energy desired energy estimate (2.19).

Even though the above PDE energy law is straightforward, the nonlinear terms in $\mu$ involves the second order derivatives, and it is not convenient to use them as test functions in numerical approximations, making it difficult to prove the discrete energy dissipation law. To overcome this difficulty, we have to reformulate the momentum equation (2.1) in an alternative form which is convenient for numerical approximation, we let $\dot{\phi} = \phi_t + \mathbf{u} \cdot \nabla\phi$, and notice that $\mu = \frac{1}{-M}\dot{\phi}$, then (2.1) can be rewritten as the following equivalent form,

$$\mathbf{u}_t + (\mathbf{u} \cdot \nabla)\mathbf{u} = \nu\Delta\mathbf{u} - \nabla p - \frac{\lambda}{M}\dot{\phi}\nabla\phi. \tag{2.27}$$

This equivalent form (2.27)-(2.2)-(2.3)-(2.4) and (2.13), (2.15), (2.16) still admits the similar energy law. Taking the inner product of (2.27) with $\mathbf{u}$, (2.13) with $\frac{\lambda}{M}\phi_t$ and (2.16) with $\lambda\phi_t$, we derive

$$\frac{1}{2}\frac{d}{dt}\|\mathbf{u}\|^2 = \nu(\partial_{\mathbf{n}}\mathbf{u}, \mathbf{u})_{\partial\Omega} - \nu\|\nabla\mathbf{u}\|^2 - \frac{\lambda}{M}(\dot{\phi}\nabla\phi, \mathbf{u}). \tag{2.28}$$

$$\frac{\lambda}{M}\|\dot{\phi}\|^2 - \frac{\lambda}{M}(\dot{\phi}, \mathbf{u} \cdot \nabla\phi) = -\lambda(\mu, \phi_t). \tag{2.29}$$

$$\lambda(\mu, \phi_t) = -\lambda\varepsilon(\partial_{\mathbf{n}}\phi, \phi_t)_{\partial\Omega} + \frac{1}{2}\lambda\varepsilon\frac{d}{dt}\|\nabla\phi\|^2 + \lambda\frac{d}{dt}(F(\phi), 1). \tag{2.30}$$

Taking the summation of the above equalities, we have

$$\begin{aligned}
\frac{1}{2}\frac{d}{dt}\|\mathbf{u}\|^2 + \frac{1}{2}\lambda\varepsilon\frac{d}{dt}\|\nabla\phi\|^2 + \lambda\frac{d}{dt}(F(\phi), 1) &= -\nu\|\nabla\mathbf{u}\|^2 - \frac{\lambda}{M}\|\dot{\phi}\|^2 \\
&\quad + \nu(\partial_{\mathbf{n}}\mathbf{u}, \mathbf{u})_{\partial\Omega} + \lambda\varepsilon(\partial_{\mathbf{n}}\phi, \phi_t)_{\partial\Omega}. \tag{2.31}
\end{aligned}$$

Using (2.25) and (2.26), we have the energy dissipation law.

$$\frac{dE_{tot}}{dt} = -\nu\|\nabla\mathbf{u}\|^2 - \frac{\lambda}{M}\|\dot{\phi}\|^2 - \lambda\gamma\|L(\phi)\|_{\partial\Omega}^2 - \|l(\phi)^{\frac{1}{2}}\mathbf{u}_s\|_{\partial\Omega}^2 - (l(\phi)\mathbf{u}_s, \mathbf{u}_w)_{\partial\Omega}. \tag{2.32}$$

We emphasize that the above derivation is suitable in a finite dimensional approximation since the test function $\phi_t$ is in the same subspaces as $\phi$. Hence, it allows us to design numerical schemes which satisfy the energy dissipation law in the discrete level.

In Shen et al. paper [41] which focus on the NSCH system, for the SCLC (2.11) and NBC (2.12), a linear *decoupled* energy stable scheme is developed. Meanwhile, for the DCLC (2.8) and GNBC (2.4), in order to obtain the energy stability, only a linear *coupled* scheme can be developed. To the best of the author's knowledge, there exist some essential difficulties to obtain the estimates for a boundary term, thus it seems formidable to obtain the decoupled energy stable scheme for DCLC coupled with GNBC.

Similar to [41], in this paper, we only develop the numerical scheme for the NSAC system with SCLC (2.11) and NBC (2.12). For the convenience of reading, we list the whole model as follows.

The hydrodynamics equation with NBC

$$\mathbf{u}_t + (\mathbf{u}\cdot\nabla)\mathbf{u} = \nu\Delta\mathbf{u} - \nabla p - \frac{\lambda}{M}\dot{\phi}\nabla\phi, \tag{2.33}$$

$$\nabla\cdot\mathbf{u} = 0, \tag{2.34}$$

$$\mathbf{u}\cdot\mathbf{n} = 0, \quad \text{on} \quad \partial\Omega. \tag{2.35}$$

$$l(\phi)(\mathbf{u}_\tau - \mathbf{u}_w) + \nu\frac{\partial\mathbf{u}_\tau}{\partial\mathbf{n}} = 0, \quad \text{on} \quad \partial\Omega. \tag{2.36}$$

The Allen-Cahn type phase field equation with SCLC

$$\phi_t + (\mathbf{u}\cdot\nabla)\phi = M(\varepsilon\Delta\phi - f(\phi) - \xi(t)), \tag{2.37}$$

$$\frac{d}{dt}\int_\Omega \phi dx = 0, \tag{2.38}$$

$$L(\phi) = 0, \quad \text{on} \quad \partial\Omega. \tag{2.39}$$

**Remark 2.2.** *There exists some essential difficulties to obtain the estimates for a boundary term for the GNBC and DCLC. In order to obtain the decoupled energy stable scheme for the system, we only consider the NBC and SCLC here.*

## 3. A Linear, Decoupled Energy Stable Numerical Scheme

We aim to construct a linear, decoupled energy stable scheme. We know that $g(\phi)$ is bounded and give the following definition

$$\bar{L} := \max_{\phi\in\mathbb{R}}|g''(\phi)| = \frac{\sqrt{2}\pi^2}{12}|\cos\theta_s|. \tag{3.1}$$

In this section we want to construct an energy stable linear scheme which is first-order temporal accurate. One introduces a function $q = (\phi^2-1)/\varepsilon$ such that one can write $f(\phi) = q\phi$. Then it follows that $\partial_t q = 2\phi\partial_t\phi/\varepsilon$. By using the function $q$, the total energy can be written as

$$E_{tot} = \frac{1}{2}\|\mathbf{u}\|^2 + \lambda\varepsilon\frac{\|\nabla\phi\|^2}{2} + \lambda\varepsilon\frac{\|q\|^2}{4} + \lambda(g(\phi), 1)_{\partial\Omega}. \tag{3.2}$$

The energy laws (2.19) and (2.32) still hold.

Let $\delta t > 0$ be a time discretization step and suppose $\mathbf{u}^n$, $\phi^n$ and $p^n$ are given, where superscript $n$ on variables denotes approximations of corresponding variables at time $n\delta t$.

A first-order linear decoupled (LD) scheme for solving the PDE system (2.33)-(2.39) is constructed as follows.

Given the initial conditions $\phi^0$, $\mathbf{u}^0$, $q^0 = ((\phi^0)^2 - 1)/\varepsilon$ and $p^0 = 0$, we aim to compute $\phi^{n+1}$, $q^{n+1}$, $\tilde{\mathbf{u}}^{n+1}$, $\mathbf{u}^{n+1}$, $p^{n+1}$ by having computed $\phi^n$, $q^n$, $\mathbf{u}^n$, $p^n$ for $n \geq 0$.

Step 1:

$$\dot{\phi}^{n+1} = M(\varepsilon\Delta\phi^{n+1} - \phi^n q^{n+1} - \xi^n), \tag{3.3}$$

$$\int_\Omega (\phi^{n+1} - \phi^n)dx = 0, \tag{3.4}$$

$$\frac{\varepsilon}{2}\frac{q^{n+1} - q^n}{\delta t} = \phi^n\frac{\phi^{n+1} - \phi^n}{\delta t}, \tag{3.5}$$

with boundary condition

$$\tilde{L}^{n+1} = \varepsilon\partial_{\mathbf{n}}\phi^{n+1} + g^{'}(\phi^n) + S(\phi^{n+1} - \phi^n) = 0, \tag{3.6}$$

where

$$\mathbf{u}_\star^n = \mathbf{u}^n - \frac{\lambda}{M}\delta t\dot{\phi}^{n+1}\nabla\phi^n, \quad \dot{\phi}^{n+1} = \frac{\phi^{n+1} - \phi^n}{\delta t} + (\mathbf{u}_\star^n \cdot \nabla)\phi^n. \tag{3.7}$$

Step 2:

$$\frac{\tilde{\mathbf{u}}^{n+1} - \mathbf{u}_\star^n}{\delta t} - \nu\Delta\tilde{\mathbf{u}}^{n+1} + B(\mathbf{u}^n, \tilde{\mathbf{u}}^{n+1}) + \nabla p^n = 0, \tag{3.8}$$

with boundary condition

$$\tilde{\mathbf{u}}^{n+1} \cdot \mathbf{n} = 0, \quad \nu\partial_{\mathbf{n}}\tilde{\mathbf{u}}_\tau^{n+1} + l(\phi^n)\tilde{\mathbf{u}}_s^{n+1} = 0, \tag{3.9}$$

where

$$B(\mathbf{u}, \mathbf{v}) = (\mathbf{u} \cdot \nabla)\mathbf{v} + \frac{1}{2}(\nabla \cdot \mathbf{u})\mathbf{v}. \tag{3.10}$$

Step 3:

$$\frac{\mathbf{u}^{n+1} - \tilde{\mathbf{u}}^{n+1}}{\delta t} + \nabla(p^{n+1} - p^n) = 0, \tag{3.11}$$

$$\nabla \cdot \mathbf{u}^{n+1} = 0, \tag{3.12}$$

$$\mathbf{u}^{n+1} \cdot \mathbf{n} = 0, \quad \text{on} \quad \partial\Omega. \tag{3.13}$$

In the above, $S$ is a parameter to be determined.

Several remarks are in order.

**Remark 3.1.** *A pressure-correction scheme is used to decouple the computation of the pressure from that of the velocity. About the overview of the projection type methods, we refer to [6,7,34].*

**Remark 3.2.** *$B(\mathbf{u}, \mathbf{v})$ is the skew-symmetric form of the nonlinear advection term is the Navier-Stokes equation, which is first introduced by Témam [45]. If the velocity is divergence* free, then $B(\mathbf{u}, \mathbf{u}) = (\mathbf{u} \cdot \nabla)\mathbf{u}$. In our scheme $\tilde{\mathbf{u}}^{n+1}$ is not divergence free, but we notice that the following identity

$$(B(\mathbf{u}, \mathbf{v}), \mathbf{v}) = 0, \quad if \quad \mathbf{u} \cdot \mathbf{n}|_{\partial\Omega} = 0. \tag{3.14}$$

Thus this identity holds regardless of whether $\mathbf{u}$ or $\mathbf{v}$ are divergence free or not, which would help to preserve the discrete energy stability.

**Remark 3.3.** Inspired by [3, 4, 26, 38, 39, 41], we introduce the explicit convective velocity $\mathbf{u}_\star^n$ in (3.3) by combining the term $\mathbf{u}^n$ with surface tension term $\dot{\phi}\nabla\phi$. From (3.7), we obtain

$$\mathbf{u}_\star^n = \tilde{B}^{-1}(\mathbf{u}^n - \frac{\phi^{n+1} - \phi^n}{M/\lambda}\nabla\phi^n), \tag{3.15}$$

where $\tilde{B} = (I + \frac{\delta t}{M/\lambda}\nabla\phi^n\nabla\phi^n)$. It is easy to get the $det(I + c\nabla\phi\nabla\phi) = 1 + c\nabla\phi \cdot \nabla\phi$, thus $\tilde{B}$ is invertible.

**Remark 3.4.** The scheme (3.3)-(3.13) is totally decoupled, linear, and first-order in time. (3.11) can be reformulated as a Poisson equation for $p^{n+1} - p^n$. Therefor, at each time step, one only needs to solve a sequence of a decoupled elliptic equations which can be solved very efficiently.

Next we shall show below, the above scheme is unconditionally energy stable.

**Theorem 3.1.** Assuming $\mathbf{u}_w = \mathbf{0}$ and $S \geq \bar{L}/2$, the solution of (3.3)-(3.13) satisfies the following discrete energy law

$$\begin{aligned}
E_{tot}^{n+1} + \frac{\delta t^2}{2}\|\nabla p^{n+1}\|^2 &+ \delta t\Big[\nu\|\nabla\tilde{\mathbf{u}}^{n+1}\|^2 + \frac{\lambda}{M}\|\dot{\phi}^{n+1}\|^2 + \|l^{1/2}(\phi^n)\tilde{\mathbf{u}}_s^{n+1}\|_{\partial\Omega}^2\Big] \\
&\leq E_{tot}^n + \frac{\delta t^2}{2}\|\nabla p^n\|^2, \quad n = 0, 1, 2, \cdots,
\end{aligned} \tag{3.16}$$

where $E_{tot}^n = \|\mathbf{u}^n\|^2/2 + \lambda(\varepsilon\|\nabla\phi^n\|^2/2 + \varepsilon\|q^n\|^2/4 + (g(\phi^n, 1))_{\partial\Omega})$.

*Proof.* By taking the $L^2$ inner product of (3.3) with $\frac{\lambda}{M}\frac{\phi^{n+1} - \phi^n}{\delta t}$, and performing integration by parts, we obtain

$$\begin{aligned}
\frac{\lambda}{M}\|\dot{\phi}^{n+1}\|^2 &- \frac{\lambda}{M}(\dot{\phi}^{n+1}, \mathbf{u}_\star^n \cdot \nabla\phi^n) - \frac{\lambda}{\delta t}(\varepsilon\partial_\mathbf{n}\phi^{n+1}, \phi^{n+1} - \phi^n)_{\partial\Omega} \\
&+ \frac{\lambda\varepsilon}{\delta t}(\frac{1}{2}\|\nabla\phi^{n+1}\|^2 - \frac{1}{2}\|\nabla\phi^n\|^2 + \|\nabla(\phi^{n+1} - \phi^n)\|^2) + \frac{\lambda}{\delta t}(\phi^n q^{n+1}, \phi^{n+1} - \phi^n) = 0,
\end{aligned} \tag{3.17}$$

where we have used the identity

$$(a - b, 2a) = |a|^2 - |b|^2 + |a - b|^2. \tag{3.18}$$

The boundary term

$$\begin{aligned}
-\frac{\lambda}{\delta t}(\varepsilon\partial_\mathbf{n}\phi^{n+1}, \phi^{n+1} - \phi^n)_{\partial\Omega} = \frac{\lambda}{\delta t}[&(g(\phi^{n+1}), 1)_{\partial\Omega} + (g(\phi^n), 1)_{\partial\Omega} \\
&- (S - \frac{g^{''}(\zeta)}{2}, (\phi_{n+1} - \phi_n)^2)_{\partial\Omega}],
\end{aligned} \tag{3.19}$$

where we have used the Taylor-expansion

$$g^{'}(\phi^n)(\phi^{n+1} - \phi^n) = g(\phi^{n+1}) - g(\phi^n) - \frac{g^{''}(\zeta)}{2}(\phi^{n+1} - \phi^n)^2. \tag{3.20}$$

Taking the $L^2$ inner product of (3.5) with $\lambda q^{n+1}$, one obtains

$$\frac{\lambda \varepsilon}{4 \delta t}(||q^{n+1}||^2 - ||q^n||^2 + ||q^{n+1} - q^n||^2) - \frac{\lambda}{\delta t}(\phi^n(\phi^{n+1} - \phi^n), q^{n+1}) = 0. \tag{3.21}$$

By taking the inner product of equation (3.7) with $\mathbf{u}_\star^n/\delta t$, we obtain

$$\frac{1}{2\delta t}(||\mathbf{u}_\star^n||^2 - ||\mathbf{u}^n||^2 + ||\mathbf{u}_\star^n - \mathbf{u}^n||^2) = -\frac{\lambda}{M}(\dot{\phi}^{n+1} \cdot \nabla \phi^n, \mathbf{u}_\star^n). \tag{3.22}$$

Taking the $L^2$ inner product of (3.8) with $\tilde{\mathbf{u}}^{n+1}$, performing integration by parts, and using identity (3.14), we obtain

$$\frac{1}{2\delta t}(||\tilde{\mathbf{u}}^{n+1}||^2 - ||\mathbf{u}_\star^n||^2 + ||\tilde{\mathbf{u}}^{n+1} - \mathbf{u}_\star^n||^2) + \nu||\nabla \tilde{\mathbf{u}}||^2 + (\nabla p^n, \tilde{\mathbf{u}}^{n+1})$$
$$- \nu(\partial_{\mathbf{n}} \tilde{\mathbf{u}}^{n+1}, \tilde{\mathbf{u}}^{n+1})_{\partial\Omega} = 0. \tag{3.23}$$

For the boundary term in the above equation, noticing that $\tilde{\mathbf{u}}_\tau^{n+1} - \tilde{\mathbf{u}}_w^{n+1} = \tilde{\mathbf{u}}_s^{n+1}$, $\tilde{\mathbf{u}}_w^{n+1} = 0$, we have

$$\nu(\partial_{\mathbf{n}} \tilde{\mathbf{u}}^{n+1}, \tilde{\mathbf{u}}^{n+1})_{\partial\Omega} = -||l^{1/2}(\phi^n)\tilde{\mathbf{u}}_s^{n+1}||^2_{\partial\Omega}. \tag{3.24}$$

By taking the $L^2$ inner product of (3.11) with $\mathbf{u}^{n+1}$ and performing integration by parts, and using the incompressible condition (3.12), we have

$$\frac{1}{2\delta t}(||\mathbf{u}^{n+1}||^2 - ||\tilde{\mathbf{u}}^{n+1}||^2 + ||\mathbf{u}^{n+1} - \tilde{\mathbf{u}}^{n+1}||^2) = 0 \tag{3.25}$$

By taking the inner product of (3.11) with $\delta t \nabla p^n$, and using the condition (3.12), we obtain

$$\frac{\delta t}{2}(||\nabla p^{n+1}||^2 - ||\nabla p^n||^2 - ||\nabla(p^{n+1} - p^n)||^2) - (\tilde{\mathbf{u}}^{n+1}, \nabla p^n) = 0. \tag{3.26}$$

We can also obtain the following equation directly from the equation (3.11),

$$\frac{\delta t}{2}||\nabla p^{n+1} - \nabla p^n||^2 = \frac{1}{2\delta t}||\tilde{\mathbf{u}}^{n+1} - \mathbf{u}^{n+1}||^2. \tag{3.27}$$

Hence, combing (3.17), (3.19), (3.21), (3.22), (3.23), (3.24), (3.25), (3.26) and (3.27), we get

$$\frac{1}{\delta t}(E_{tot}^{n+1} - E_{tot}^n) + \frac{\delta t}{2}(||\nabla p^{n+1}||^2 - ||\nabla p^n||^2) =$$
$$- \left[\nu||\nabla \tilde{\mathbf{u}}^{n+1}||^2 + \frac{\lambda}{M}||\dot{\phi}^{n+1}||^2 + ||l^{1/2}(\phi^n)\tilde{\mathbf{u}}_s^{n+1}||^2_{\partial\Omega}\right]$$
$$- \frac{1}{2\delta t}(||\tilde{\mathbf{u}}^{n+1} - \mathbf{u}_\star^n||^2 + ||\mathbf{u}_\star^n - \mathbf{u}^n||^2) - \frac{\lambda \varepsilon}{4\delta t}||q^{n+1} - q^n||^2.$$
$$- \frac{\lambda}{\delta t}(S - \frac{g^{''}(\zeta)}{2}, (\phi^{n+1} - \phi^n)^2)_{\partial\Omega}. \tag{3.28}$$

Thus, by the assumption $S \geq \bar{L}/2$ we get the desired energy law (3.16).

# 4. Spatial discretization

Since the proofs of energy stability for scheme is based on weak form and integration by parts, a suitable spatial discretization method can be adopted to keep the energy dissipation properties of the semi-discretized scheme.

In this section, we consider the domain $\Omega = [0, L_x] \times [0, L_y]$, the boundary $\partial\Omega = \Gamma_1 \cup \Gamma_2 \cup \Gamma_3 \cup \Gamma_4$, and the definition of the boundary is showed in Fig. 4.1. We take the implementation of a finite element method for the spatial discretization to test the approximation properties of our time discretization LD scheme.

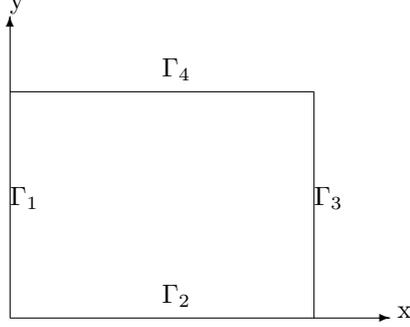

Fig. 4.1. A diagram of the boundary.

## 4.1. Fully discrete schemes and energy stability

Now we consider the fully discrete scheme to solve the PDE system. Let $S_h \subset H^1(\Omega)$ be the finite-dimensional subspace, which is constructed by the piecewise linear functions.

We define $\tilde{S}_h^0 = \{u : u|_{\Gamma_1 \cup \Gamma_3} = 0, u \in S_h\}$, $S_h^0 \subset H_0^1(\Omega)$. Let $V_{\mathbf{u}_h} = \tilde{S}_h^0 \times S_h^0$ and $M_h \subset L_0^2(\Omega)$ be two finite-dimensional spaces satisfying the inf-sup condition,

$$\inf_{q_h \in M_h} \sup_{\mathbf{u}_h \in V_{\mathbf{u}_h}} \frac{\int_\Omega q_h \nabla \cdot \mathbf{u}_h dx}{\|q_h\| \|\mathbf{u}_h\|_1} \geq C, \tag{4.1}$$

where $C > 0$ independent of mesh size $h$ and $\|\mathbf{u}_h\|_1 = \|\nabla \mathbf{u}_h\| + \|\mathbf{u}_h\|$.

In the framework of the finite element spaces, the scheme (3.3)-(3.13) reads as follows.

Find $(\phi_h^{n+1}, q_h^{n+1}, \tilde{\mathbf{u}}_h^n, p_h^{n+1}, \mathbf{u}_h^{n+1}) \in S_h \times S_h \times V_{\mathbf{u}_h} \times M_h \times V_{\mathbf{u}_h}$, such that for all $(\psi_h, \mathbf{v}_h, g_h) \in S_h \times V_{\mathbf{u}_h} \times M_h$ there holds

$$\begin{aligned}
(\frac{\phi^{n+1} - \phi^n}{\delta t}, \psi_h) + (\tilde{\mathbf{u}}_{\star h}^n \cdot \nabla \phi_h^n, \psi_h) + M\varepsilon(\nabla \phi_h^{n+1}, \nabla \psi_h) - M(\varepsilon \partial_{\mathbf{n}} \phi_h^{n+1}, \psi_h)_{\partial\Omega} \\
+ M(\phi_h^n [q_h^n + \frac{2}{\varepsilon} \phi_h^n (\phi_h^{n+1} - \phi_h^n)], \psi_h) + M(\xi^n, \psi_h) = 0,
\end{aligned} \tag{4.2}$$

$$\begin{aligned}
(\frac{\tilde{\mathbf{u}}_h^{n+1} - \mathbf{u}_{\star h}^n}{\delta t}, \mathbf{v}_h) + (B(\mathbf{u}_h^n, \tilde{\mathbf{u}}_h^{n+1}), \mathbf{v}_h) + \nu(\nabla \tilde{\mathbf{u}}_h^{n+1}, \nabla \mathbf{v}_h) \\
- \nu(\partial_{\mathbf{n}} \tilde{\mathbf{u}}_h^{n+1}, \mathbf{v}_h)_{\partial\Omega} + (\nabla p_h^n, \mathbf{v}_h)) = 0,
\end{aligned} \tag{4.3}$$

$$(\frac{\mathbf{u}_h^{n+1} - \tilde{\mathbf{u}}_h^{n+1}}{\delta t}, \mathbf{v_h}) + (\nabla(p_h^{n+1} - p_h^n), \mathbf{v}_h) = 0, \tag{4.4}$$

$$(\nabla \cdot \mathbf{u}_h^{n+1}, g_h) = 0, \tag{4.5}$$

$$\frac{\varepsilon}{2}(\frac{q_h^{n+1} - q_h^n}{\delta t}, \psi_h) = (\phi_h^n \frac{\phi_h^{n+1} - \phi_h^n}{\delta t}, \psi_h), \tag{4.6}$$

where

$$-\varepsilon\partial_{\mathbf{n}}\phi_h^{n+1} = g^{'}(\phi_h^n) + S(\phi_h^{n+1} - \phi_h^n), \tag{4.7}$$

$$-\nu\partial_{\mathbf{n}}\tilde{\mathbf{u}}_h^{n+1} = l(\phi_h^n)\tilde{\mathbf{u}}_{sh}^{n+1}. \tag{4.8}$$

**Remark 4.1.** *Note that the update of $q_h^{n+1}$ in (4.6) is decoupled from the rest of equations.*

For the convenience of establishing the stability of the fully discrete scheme (4.2)-(4.6), we introduce the discrete divergence operator $B_h : V_{\mathbf{u}_h} \to M_h$ such that for $\mathbf{u}_h \in V_{\mathbf{u}_h}$ and $p_h \in M_h$

$$(B_h\mathbf{u}_h, p_h) := -(\nabla \cdot \mathbf{u}_h, p_h) = (\mathbf{u}_h, \nabla p_h) := (\mathbf{u}_h, B_h^T p_h), \tag{4.9}$$

where $B_h^T : M_h \to V_{\mathbf{u}_h}$ is the transpose of $B_h$ (the discrete gradient operator). Hence one can write the projection step (4.4) and (4.5) in discrete form

$$\mathbf{u}_h^{n+1} - \tilde{\mathbf{u}}_h^{n+1} + B_h^T(p_h^{n+1} - p_h^n) = 0, \quad \text{in} \quad M_h, \tag{4.10}$$

$$B_h\mathbf{u}_h^{n+1} = 0, \quad \text{in} \quad M_h. \tag{4.11}$$

Thus, the fully discrete energy of the system also satisfies the following law.

**Theorem 4.1.** *Given that $q_h^n \in S_h$, $\phi_h^n \in S_h$, $\mathbf{u}_h^n \in V_{\mathbf{u}_h}$, and $p_h^n \in M_h$, the system (4.2)-(4.6) admits a unique solution $(\phi_h^{n+1}, q_h^{n+1}, \mathbf{u}_h^{n+1}, p_h^{n+1}) \in S_h \times S_h \times V_{\mathbf{u}_h} \times M_h$ at the time $t^{n+1}$ for any parameters $h > 0$ and $\delta t > 0$. Moreover, by the assumption of $\mathbf{u}_w = 0$ and $S \geq \bar{L}/2$, the solution satisfies a discrete energy law*

$$E_{tot}_h^{n+1} + \frac{\delta t^2}{2}\|B_h^T p_h^{n+1}\|^2 + \delta t\Big[\nu\|\nabla\tilde{\mathbf{u}}_h^{n+1}\|^2 + \frac{\lambda}{M}\|\dot{\phi}_h^{n+1}\|^2 + \|l^{1/2}(\phi_h^n)\tilde{\mathbf{u}}_{sh}^{n+1}\|_{\partial\Omega}^2\Big]$$

$$\leq E_{tot}_h^n + \frac{\delta t^2}{2}\|B_h^T p_h^n\|^2, \quad n = 0, 1, 2, \cdots, \tag{4.12}$$

*where $E_{tot}_h^n = \|\mathbf{u}_h^n\|^2/2 + \lambda(\varepsilon\|\nabla\phi_h^n\|^2/2 + \varepsilon\|q_h^n\|^2/4 + (g(\phi_h^n), 1))_{\partial\Omega})$, and $\dot{\phi}_h^{n+1} = (\phi_h^{n+1} - \phi_h^n)/\delta t + (\mathbf{u}_{\star h}^{n+1} \cdot \nabla)\phi_h^n$.*

*Proof.* Noting that the scheme (4.2)-(4.6) is a linear system, thus the unique solvability would follow from the energy law (4.12).

By working with (4.2), (4.3), (4.6), (4.10) and (4.11), and by the same proof of Theorem 3.1, we can easily prove that it satisfies the energy law (4.12).

## 5. Numerical simulation

We present in this section some numerical experiments using the schemes constructed in the section 3.

Our spatial discretization is based on the finite element method. We use the inf-sup stable $Iso-P2/P1$ element [44] for the velocity and pressure, and linear element for the phase function $\phi$.

**Example 1: Accuracy test.** We first test the convergence rates of the proposed schemes (3.3)-(3.13). In $\Omega = [0, 2]^2$, we set the exact solution as

$$\begin{cases} \phi(t, x, y) = 2 + \cos(\pi x)\cos(\pi y)\sin t, \\ u(t, x, y) = \pi\sin(2\pi y)\sin^2(\pi x)\sin t, \\ v(t, x, y) = -\pi\sin(2\pi x)\sin^2(\pi y)\sin t, \\ p(t, x, y) = \cos(\pi x)\sin(\pi y)\sin t. \end{cases} \tag{5.1}$$

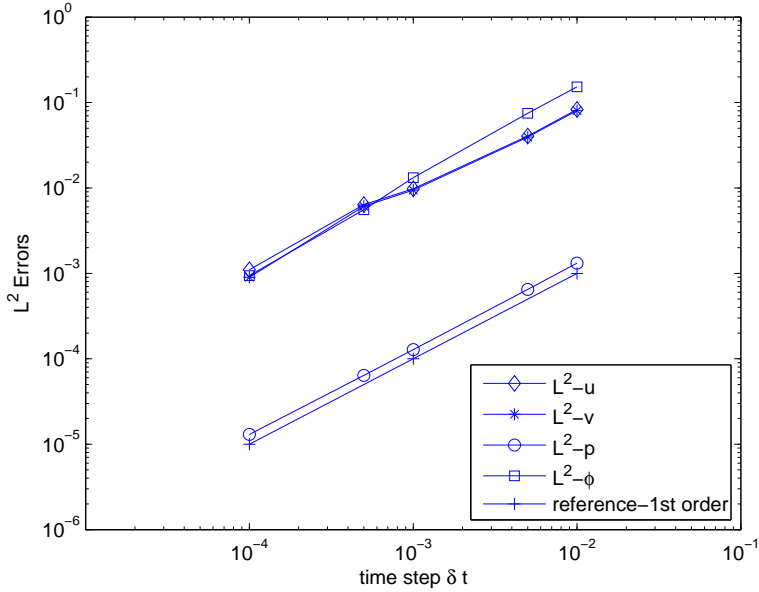

Fig. 5.1. **Example 1**. Temporal convergence rates of $L^2$ errors for the velocity ($\mathbf{u} = (u, v)$), pressure $p$, and phase field function $\phi$ in scheme (3.3)-(3.13).

We choose $\varepsilon = 0.025$, $\nu = 1$, $M = 0.001$, $\lambda = 10^{-7}$, $\gamma = 1000$, $l(\phi) = 1/(0.19)$, $\theta_s = 90°$. Some suitable force fields are imposed such that the given solutions satisfy the coupled systems. We used 10145 nodes and 19968 triangle elements in this test.

For the accuracy test, we plot the $L^2$ errors of the velocity, pressure and phase function between the numerical solution and the exact solution at $t = 1$ with different time step size $\delta t = 0.0001$, $0.0005$, $0.001$, $0.005$, and $0.01$ in Figure 5.1. The numerical results show that scheme (3.3)-(3.13) is first-order accurate in time for all variables.

**Example 2: Immiscible Couette flow.**

In this case the domain is set to be $\Omega = [0, 4] \times [0, 0.8]$ and the solid wall velocity is set to be $V = 0.7$ and $\lambda = 0.1$. We choose the different value of $\theta_s$ and the other parameters in **Example 1**. We set the initial values $\mathbf{u}^0 = 0$ and

$$\phi^0(x, y) = \begin{cases} 1, & (x, y) \in \Gamma, \\ -1, & \text{else}, \end{cases} \tag{5.2}$$

where $\Gamma$ is the initial configuration of fluid 1, and $\Gamma = \{(x, y) \in \Omega | |x - 2| \leq 1\}$. We obtain 9153 nodes and 17920 triangle elements in the fixed domain produced by the PDE tool in Matlab. The Couette flow is generated by moving the top and bottom walls at a speed $V$ along $\pm x$, respectively.

In **Example 2 (a)**, we set the contact angle $\theta_s = 60°$, and we see that the contact line of Fig. 5.2 is moving from $T = 0.01$ to $T = 3$. The contact line goes to steady state form $T = 3$ to $T = 5$.

In **Example 2 (b)**, we set the angle $\theta_s = 103°$ on the top boundary and $\theta_s = 77°$ on the bottom boundary. In Fig. 5.3, we also see that the contact line is moving from $T = 0.01$ to $T = 3$ and reaches equilibrium state at $T = 5$. Theses results is consistent to the numerical results in [24, 41].

**Example 3: Spreading and dewetting of a droplet.**

In this example we set the domain $\Omega = [0,4] \times [0,1.2]$ and the initial configuration $\Gamma = \{(x,y) \in \Omega | (x-2)^2 + y^2 \leq 0.64\}$.

We want to simulate the process of the spreading and dewetting of a droplet. The solid wall velocity is set to be $V = 0$ and other parameters are set to be the same in **Example 2**. In this example we select the contact angle $\theta_s = 30°$ in the wetting case and $\theta_s = 150°$ in the dewetting case. We see that the contact line is spreading out from $T = 0.01$ to $T = 0.5$ in the right of Fig. 5.4 and goes to steady state at $T = 5$ in the right of Fig. 5.5. In the dewetting case the contact line is shrinking from $T = 0.01$ to $T = 0.5$ in the right of Fig. 5.4 and also goes to steady state at $T = 5$ in the right of Fig. 5.5. Fig. 5.6 shows that the energy of the two cases is decreasing and goes to steady state at the end. Due to the boundary energy $E_s(\phi) = \lambda(g(\phi), 1)_{\partial\Omega}$, the total energy in dewetting case is less than zero and that in wetting case is bigger than zero. Both cases match to the numerical results in [24,41].

**Example 4: A droplet in a cylinder with a hole.**

We now consider the physical domain $\Omega$ to be a cylinder with a hole, and assume that all variables are axisymmetric. For the implementation of the cylindrical coordinates $(r, \theta, z)$, the Navier-Stokes equations in the system read as in [23,32]. Without considering the azimuthal $\theta$ direction, the system can be considered to be two-dimensional problems briefly.

In this case we select the contact angle $\theta_s = 30°$ for $(a)$ and $\theta_s = 150°$ for $(b)$, considering the effect of gravity on the system. We can use the momentum equation as follows,

$$\mathbf{u}_t + (\mathbf{u} \cdot \nabla)\mathbf{u} + \nabla p = \lambda\mu\nabla\phi + \phi g_0 e_z, \tag{5.3}$$

where $g_0$ is the gravity acceleration and $e_z = (0,1)$. The computation domain is $R = 1$, $Z = 2$, which has a rectangle hole. The irregular domain is discretized into 20105 nodes and 39680 triangle elements by MATLAB. The initial configuration is set to be a hole around a semicircle.

The Fig. 5.7 and Fig. 5.8 show that the droplet is dropping with time going on and the contact line is spreading. Similarly, the Fig. 5.9 and Fig. 5.10 show that the droplet is dropping and the contact line is shrinking.

## 6. Conclusions and remarks

In this paper we studied the Navier-Stokes coupled with Allen-Cahn phase field model for two phase immiscible complex fluids. With the implementation of the general Navier boundary condition and the dynamic condition, we derived the energy dissipation law for the model. We then constructed a linear decoupled energy stable scheme by the introduction of the auxiliary function.

We have not only considered time discretizations here, but also take the implementation on the spatial dicretizations by the finite element method. Moreover, we prove that the fully discrete scheme is also energy stable. At last we presented ample numerical results to validate the numerical scheme and illustrate the contact lines' movements, wetting and dewetting phenomena.

## Acknowledgments

R. Chen is partially supported by the China Postdoctoral Science Foundation grant No. 2016M591122. X. Yang is partially supported by NSF DMS-1200487, NSF DMS-1418898,

AFOSR FA9550-12-1-0178. H. Zhang is partially supported by NSFC/RGC Joint Research Scheme No. 11261160486, NSFC grant No. 11471046, 11571045 and the Ministry of Education Program for New Century Excellent Talents Project NCET-12-0053.

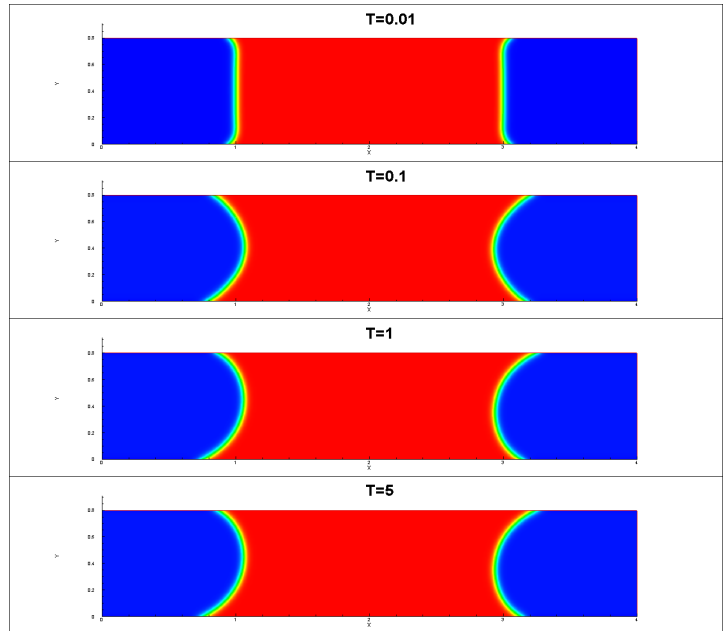

Fig. 5.2. **Example 2 (a):** The motion of the two fluids at the time of $T = 0.01$, 0.1, 1, 5 with contact angle $\theta_s = 60°$.

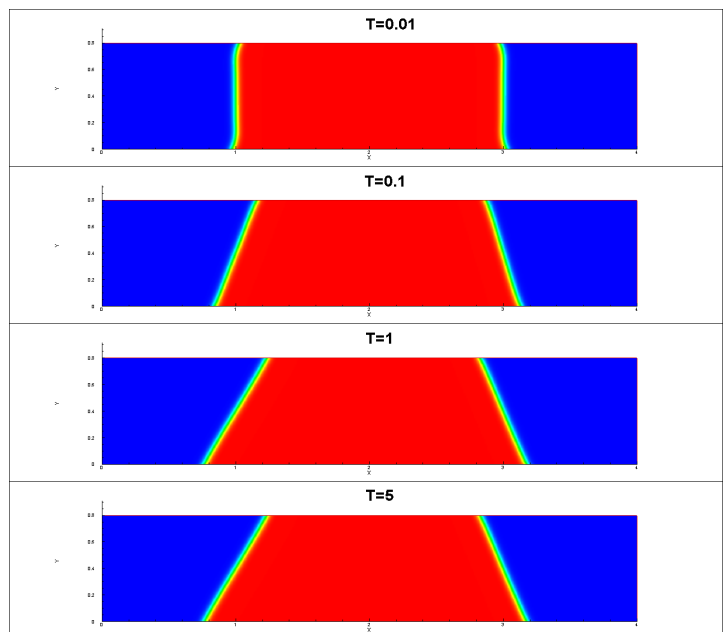

Fig. 5.3. **Example 2 (b):** The motion of the two fluids at the time of $T = 0.01$, 0.1, 1, 5 with contact angle $\theta_s = 103°$ on the top and $\theta_s = 77°$ on the bottom.

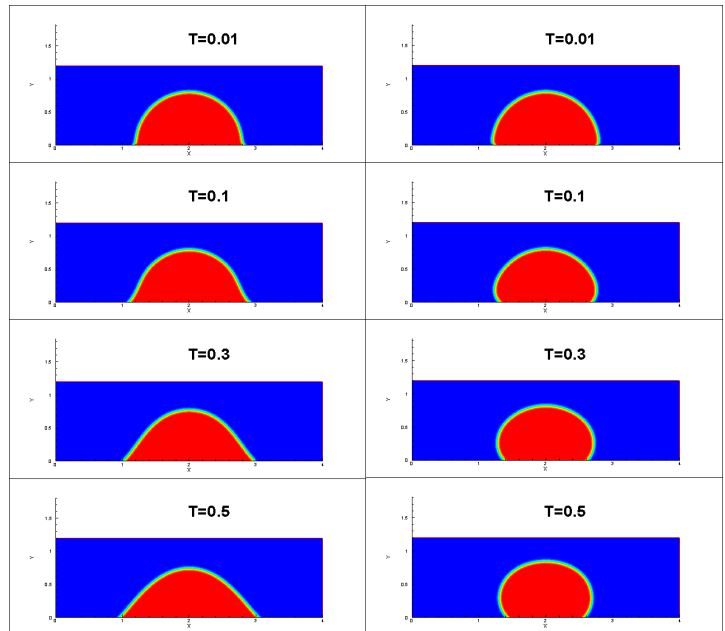

Fig. 5.4. **Example 3:** The motion of the two fluids at the time of $T = 0.01$, $0.1$, $0.3$, $0.5$ with the different contact angle. Left: $\theta_s = 30°$, right: $\theta_s = 150°$.

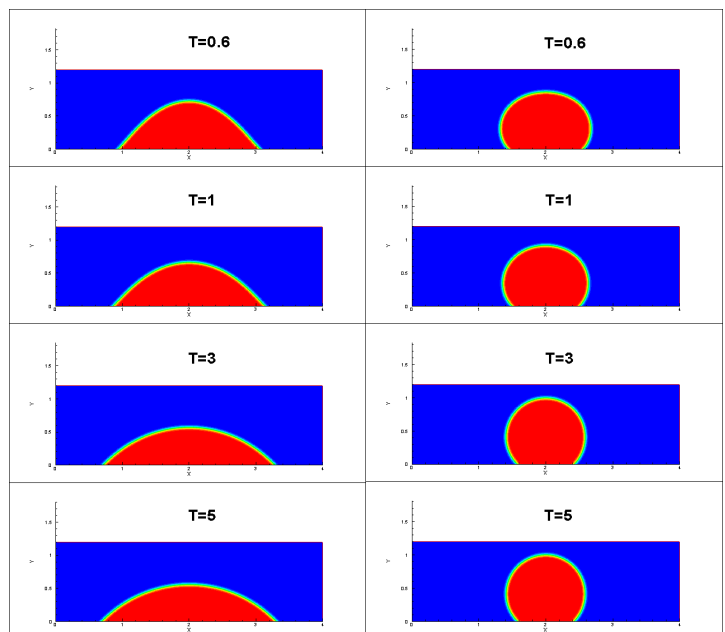

Fig. 5.5. **Example 3:** The motion of the two fluids at the time of $T = 0.6$, $1$, $3$, $5$ with the different contact angle. Left: $\theta_s = 30°$, right: $\theta_s = 150°$.

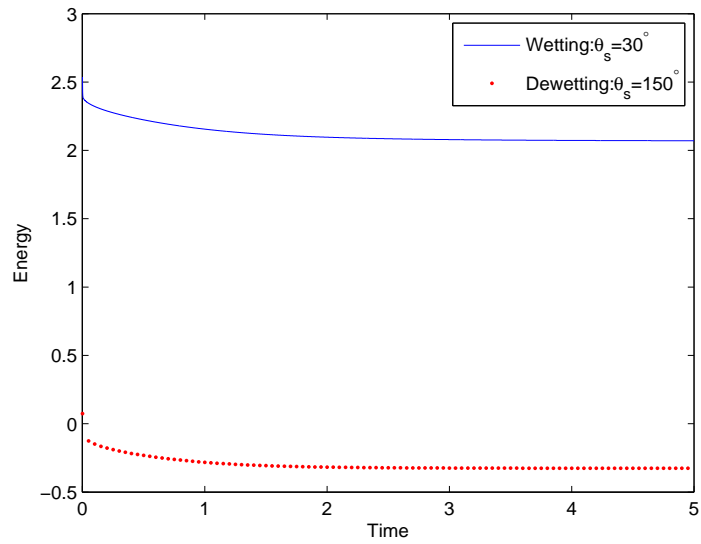

Fig. 5.6. **Example 3:** Time evolution of the free energy functional for wetting case $\theta_s = 30°$ and dewetting case $\theta_s = 150°$. The energy curves of the two cases are decreasing and reach the steady state at the end.

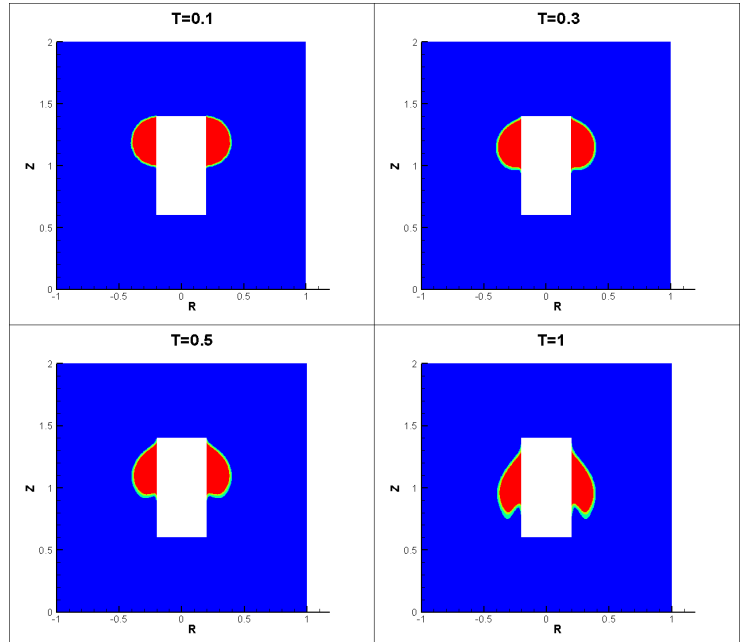

Fig. 5.7. **Example 4(a):** The motion of the two fluids at the time of $T = 0.1, 0.3, 0.5, 1$ with the contact angle $\theta_s = 30°$.

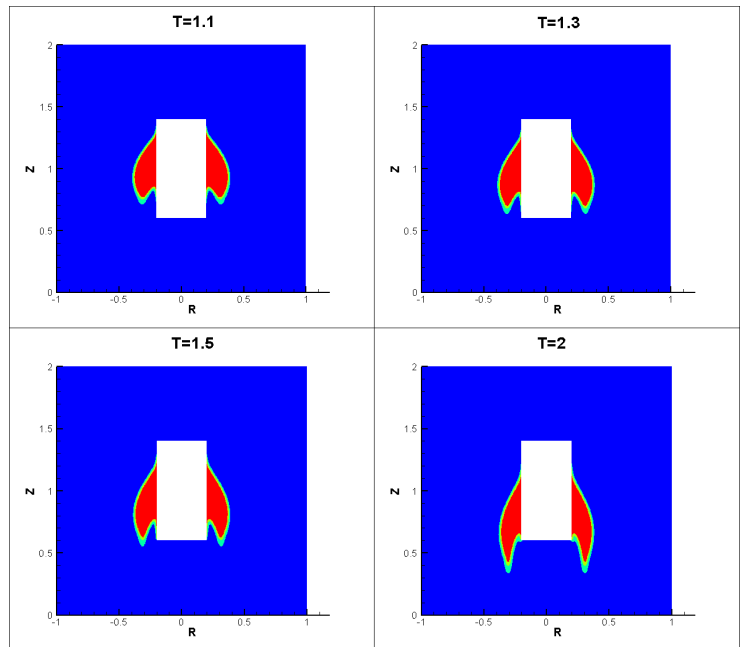

Fig. 5.8. **Example 4(a)**: The motion of the two fluids at the time of $T = 1.1, 1.3, 1.5, 2$ with the contact angle $\theta_s = 30°$.

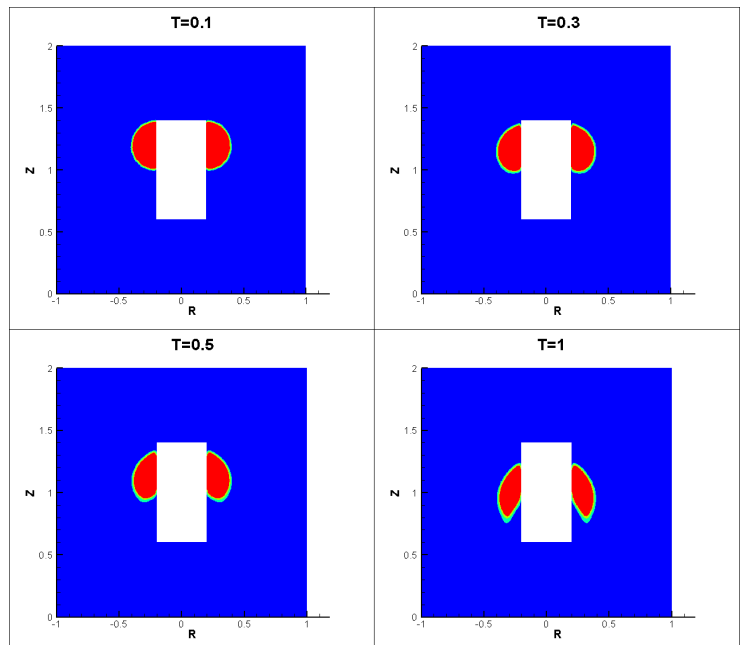

Fig. 5.9. **Example 4(b)**: The motion of the two fluids at the time of $T = 0.1, 0.3, 0.5, 1$ with the contact angle $\theta_s = 150°$.

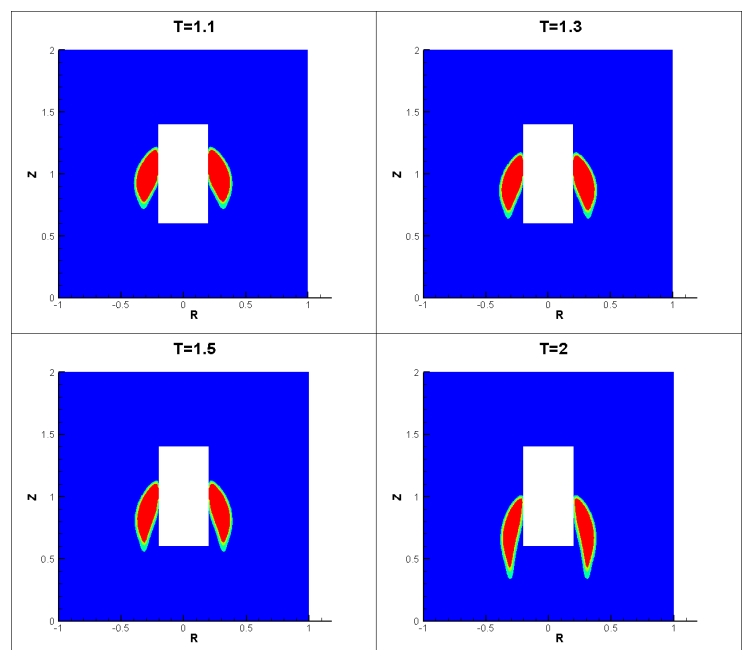

Fig. 5.10. **Example 4(b)**: The motion of the two fluids at the time of $T = 1.1, 1.3, 1.5.2$ with the contact angle $\theta_s = 150°$.